\documentclass[11pt]{article}
\usepackage{amsmath,amssymb,amscd,array}

\let\ssection=\section
\renewcommand{\section}{\setcounter{equation}{0}\ssection}

\def\d{\delta}

\def\om{\omega}
\def\r{\rho}
\def\a{\alpha}
\def\b{\beta}
\def\s{\sigma}
\def\vfi{\varphi}

\def\l{\lambda}
\def\m{\mu}

\def\implies{\Rightarrow}

\newcommand{\Vect}{\mathrm{Vect}}

\newcommand{\cqfd}{\hspace*{\fill}\rule{3mm}{3mm}}
\newcommand{\cqf}{\hspace*{\fill}\rule{2mm}{2mm}}

\begin{document}

\frenchspacing

\def\d{\delta}
\def\g{\gamma}
\def\om{\omega}
\def\r{\rho}
\def\a{\alpha}
\def\b{\beta}
\def\s{\sigma}
\def\vfi{\varphi}
\def\l{\lambda}
\def\m{\mu}
\def\implies{\Rightarrow}

\oddsidemargin .1truein
\newtheorem{thm}{Theorem}[section]
\newtheorem{lem}[thm]{Lemma}
\newtheorem{cor}[thm]{Corollary}
\newtheorem{pro}[thm]{Proposition}
\newtheorem{ex}[thm]{Example}
\newtheorem{rmk}[thm]{Remark}
\newtheorem{defi}[thm]{Definition}
\title{On $\mathfrak{sl}(2)$-equivariant quantizations}
\author{S. Bouarroudj
\thanks{E-mail:Bouarroudj.Sofiane@uaeu.ac.ae
\hfill\break U. A. E. University, Department of Mathematical
Sciences, P. O. Box 17551, Al-Ain, U. A. E.} \and M. Iyadh Ayari
\thanks{E-mail:mohdayari@uaeu.ac.ae
\hfill\break Universit\'e du centre, institut sup\'erieur
d'informatique et de gestion de Kairouan, Tunisia. } }
\date{}
\maketitle

\begin{abstract}
\noindent By computing certain cohomology of $\mathrm{Vect}(M)$ of
smooth vector fields we prove that on 1-dimensional manifolds $M$
there is no quantization map intertwining the action of
non-projective embeddings of the Lie algebra $\mathfrak{sl}(2)$ into
the Lie algebra $\mathrm{Vect} (M)$. Contrariwise, for projective
embeddings $\mathfrak{sl}(2)$-equivariant quantization exists.
\end{abstract}

%
\section{Introduction}
{\em Quantization}, according to one of many definitions, is a
procedure which to each polynomial on a symplectic manifold (in
our case it is $T^*M$) assigns a differential operator on a
Hilbert space (see, e.g., \cite{wo}). An {\em equivariant} (with
respect to a finite-dimensional Lie subalgebra
$\mathfrak{g}\subset \mathrm{Vect} (M)$ of the Lie algebra of all
smooth vector fields on $M$) quantization intertwines a local
action of ${\mathfrak g}\subset \mathrm{Vect} (M)$, see \cite{dlo,
do2, lol}. For $M=\mathbb{R}$ or $S^1$, this point of view was
initiated in \cite{cmz} in the context of a Lie subalgebra
$\mathfrak{sl}(2)\subset \mathrm{Vect} (M)$. In fact, the Lie
algebra $\mathfrak{sl}(2)$ has various embeddings into
$\mathrm{Vect}(M)$ (for $M=\mathbb{R}$, only $n=0$ is possible):
\[
\mathfrak{l}_n=\text{Span}\left (x^{-n}\frac{d}{dx},\quad
x\frac{d}{dx},\quad x^{n+2}\frac{d}{dx} \right).
\]
The Lie algebra $\mathfrak{sl}(2)$ considered in
\cite{bo, cmz, gar} corresponds to ${\mathfrak l}_0$ generated by
the natural action of $\mathrm{SL}(2)$ by linear-fractional (or
M\"{o}bius) transformations (cf. \cite{olv}). Are there any other
Lie subalgebras of $\mathrm{Vect}(M)$ for
which the equivariant quantization procedure exists? On
$\mathbb{R},$ all subalgebras that contain $\frac{d}{dx}$ were
classified by Lie \cite{li} (see also \cite{es}, \cite{post}). Apart from
${\mathfrak l}_{0}$, they are (for a fixed
$s\in\mathbb{R}\setminus \{0\}$):
$$
\begin{array}{cclccclccl}
{\mathfrak g}_{0}&=&\displaystyle \text{Span}\left
(\frac{d}{dx}\right ),& {\mathfrak h}_{0}&=&\displaystyle
\text{Span}\left (\frac{d}{d x},\quad e^{s x}\frac{d}{d x}\right
),_{1}&=&\displaystyle\text{Span}\left (\frac{d}{d
x},\quad x\frac{d}{d x}\right ),
\end{array}
$$
$$
\begin{array}{l}
{\mathfrak k}_{1}=\displaystyle \text{Span}\left (\frac{d}{d x},\;
\sin (sx)\frac{d }{d x}, \; \cos (sx)\frac{d}{d x}
\right ),\quad 
{\mathfrak k}_{2}=\displaystyle \text{Span}\left (\frac{d}{d x},
\; \sinh (sx)\frac{d}{d x},\; \cosh (sx)\frac{d}{d x}\right ).
\end{array}
$$
The Lie subalgebras ${\mathfrak l}_n,$ ${\mathfrak k}_{1}$ and
${\mathfrak k}_{2}$ correspond to various embeddings of
$\mathfrak{sl}(2).$ The vector field $\frac{d}{dx}$ spans a
commutative Lie algebra isomorphic to $\mathfrak{so}(2).$ We
require that the Lie algebras we are dealing with contain
$\mathfrak{so}(2).$ This property implies that the invariant
differential operators we are studying are with constant
coefficients (see sec. \ref{mai}).

The existence of an equivariant quantization map induces a
deformation of the $\mathrm{Vect}(M)$-action on the space of
polynomials in the sense of the Nijenhuis-Richardson deformation
theory \cite{nr}. As is well known, deformation of modules is
related to the computation of a certain first cohomology group. In
our case, the deformation becomes trivial when we restrict the
action to $\mathfrak{sl}(2);$ besides, we deal only with a
deformation that is expressed in terms of differential operators.
Therefore, we will deal with the relative cohomology group
\begin{equation}
\label{ver} {\mathrm H}^{1}(\mathrm{Vect}(M), \mathfrak{sl}(2);
\mathrm{Hom}_{\mathrm{diff}}({\cal F}_{\lambda},{\cal F}_\mu)),
\end{equation}
where ${\cal F}_{\lambda}$ stands for the space of densities of
weight $\lambda$ and $\mathrm{Hom}_{\mathrm{diff}}({\cal
F}_{\lambda},{\cal F}_\mu)$ stands for the space of differential
operators on weighted densities.

Our first main result is the computation of (\ref{ver})
for $\mathfrak{sl}(2)$ realized as $\mathfrak{k}_1$ or
$\mathfrak{k}_2.$ Here, we only deal with cochains given by smooth
maps. Using this result, we prove that there is no equivariant
quantization map intertwining the actions of the Lie subalgebras
$\mathfrak{k}_1$ and $\mathfrak{k}_2.$

Boniver and Mathonet \cite{bm} have investigated, in
multi-dimensional case, all maximal Lie subalgebras of
$\Vect(\mathbb{R}^n)$ for which the equivariant quantization
exists.
\section{Relative cohomology, statement of main results}
In what follows, $M$ stands for either $\mathbb{R}$ or $S^1$ and
$\mathfrak{k}=\mathfrak{k}_1$ or $\mathfrak{k}_2.$

The space of tensor densities, denoted by ${\cal F}_{\lambda},$ is
the space of sections of the line bundle $ (T^*M)^{\otimes
\lambda}$, where $\lambda\in \mathbb{R}$. This space has the
following  structure of a $\mathrm{Vect}(M)$-module: For any
$\varphi(dx)^\lambda \in {\cal F}_{\lambda}$ and $X\frac{d}{dx}\in
\mathrm{Vect}(M),$ we put
\begin{equation}
\label{dens} L_{X\frac{d}{dx}}^{\lambda}(\varphi (dx)^{\lambda}) =(X
{\varphi}' +\lambda \,\varphi \,\text{div}\,X)(dx)^{\lambda}.
\end{equation}
Examples of such spaces are ${\cal F}_0=C^{\infty}(M),$ ${\cal
F}_{-1}=\mathrm{Vect}(M),$ and ${\cal F}_1=\Omega^{1}(M)$.

A result of \cite{bo} is as follows:
\begin{equation} \label{mou1}
\mathrm H^1(\mathrm{Vect}(M),{\mathfrak l}_0;
\mathrm{Hom}_{\mathrm{diff}}({\cal F}_{\lambda},{\cal F}_\mu))=
\begin{cases}
\mathbb{R}& \text{if $\mu-\lambda=\begin{cases}2& \text{and $
\lambda\not =-\frac{1}{2}$}\\
3& \text{and $\lambda\not =-1$}\\
4& \text{and $\lambda\not =-\frac{3}{2}$}\end{cases}$ or } (\lambda, \mu)=\\
&(-4,1),\; (0,5), \; \left (-\frac{5\pm
\sqrt{19}}{2},
\frac{7\pm \sqrt{19}}{2}\right )\\
0& \mbox{otherwise.}
\end{cases}
\end{equation}
It is known (\cite{l}) that
\begin{equation} \label{mou2}
\mathrm H^1({\mathfrak l}_0; \mathrm{Hom}_{\mathrm{diff}}({\cal
F}_{\lambda},{\cal F}_\mu))= \left\{
\begin{array}{ll}
\mathbb{R}& \mbox{if} \quad \mu-\lambda=0\\[1mm]
\mathbb{R}^2& \mbox{if} \quad (\lambda, \mu)= \left
(\frac{1-n}{2},\frac{n+1}{2}\right ), \quad \mbox{where }
n \in\mathbb{N}\backslash \{0\}\\
0& \mbox{otherwise.}
\end{array}
\right.
\end{equation}

\begin{thm}
 \label{t1}
i) $ \mathrm{H}^{1}(\mathrm{Vect}(M),{\mathfrak k};
\mathrm{Hom}_{\mathrm{diff}}({\cal F}_{\lambda},{\cal F}_\mu))=0$
for all $\lambda$ and $\mu.$

ii) $\displaystyle \mathrm{H}^{1}({\mathfrak k};
\mathrm{Hom}_{\mathrm{diff}}({\cal F}_{\lambda},{\cal
F}_\mu))\simeq\displaystyle{\mathrm{H}^{1}(\mathrm{Vect}(M);
\mathrm{Hom}_{\mathrm{diff}}({\cal F}_{\lambda},{\cal F}_\mu))}.$
Moreover, if $M=\mathbb{R},$ we have (here
$\mathrm{H}^1_{\mathrm{diff}}$ denotes the differential cohomology
which means we deal only with cochains that are given by
differential operators)
\begin{equation}
\label{tha}
\mathrm{H}^{1}_{\mathrm{diff}}(\mathrm{Vect}(\mathbb{R}),
{\mathfrak k}; \mathrm{Hom}_{\mathrm{diff}}({\cal
F}_{\lambda},{\cal F}_\mu))= \begin{cases}
\mathbb{R}& \text{if $\begin{cases}\mu-\lambda=0, 2, 3\;\text{ or }(\lambda, \mu)=&\\
(-4,1),\; (0,5),\;\left (-\frac{5\pm \sqrt{19}}{2},
\frac{7\pm \sqrt{19}}{2}\right )&\end{cases}$}\\
\mathbb{R}^2& \text{if $(\lambda,\mu)=(0,1)$}\\
0& \text{otherwise}.
\end{cases}
\end{equation}
\end{thm}
\begin{rmk} {\rm (i) The 1-cocycles that
span $\mathrm{H}^1(\mathrm{Vect}(\mathbb{R});
\mathrm{Hom}_{\mathrm{diff}}({\cal F}_\lambda,{\cal F}_\mu))$ were
given in \cite{f}.

(ii) The cohomology group
$\mathrm{H}^{1}_{\mathrm{diff}}(\mathrm{Vect}(S^1); {\mathfrak k};
\mathrm{Hom}_{\mathrm{diff}}({\cal F}_{\lambda},{\cal F}_\mu))$
cannot be described as in (\ref{tha}). The cocycles
\[
c_1(X
\frac{d}{dx})(\phi(dx)^{\lambda})=\mathrm{div}X\,\phi(dx)^{\mu-\lambda},
\mbox{ and }\quad c_2(X \frac{d}{dx})(\phi(dx)^{\lambda})=X\cdot
\phi(dx)^{\mu-\lambda},
\]
are not trivial and not cohomologous to each other for every pair $(\lambda,\mu)$
such that $\mu-\lambda=0.$}
\end{rmk}

The proof of Theorem \ref{t1} will be the subject of subsection
\ref{sw} following the investigation of ${\mathfrak k}$-invariant
bilinear differential operators.
\section{Invariant operators and cohomology}
\subsection{\bf{${\mathfrak k}$-}invariant bilinear differential operators}
\label{mai}
\quad First, we recall Grozman's result \cite{groz,groz2} on the
classification of  bilinear differential operators on weighted
densities invariant with respect to $\Vect(M):$ \\

\noindent i) Every zero-order operator ${\cal F}_{\gamma}\otimes
{\cal F}_{\lambda} \longrightarrow {\cal F}_{\gamma+\lambda}$ is
of the form
$$
\varphi (dx)^\gamma\otimes \psi(dx)^\lambda \mapsto \varphi\cdot
\psi\,(dx)^{\gamma+\lambda}.
$$
ii) Every first-order operator ${\cal F}_{\gamma}\otimes {\cal
F}_{\lambda} \longrightarrow {\cal F}_{\gamma+\lambda+1}$ is given
by the Poisson bracket:
$$
\varphi\, (dx)^{\gamma} \otimes \psi\,(dx)^{\lambda} \mapsto\left
\{\varphi,\psi\right \}(dx)^{\gamma+\lambda+1}\quad \mbox{ where }
\quad \left \{\varphi,\psi\right \}= \gamma\, \varphi\,
\psi'-\lambda\, \varphi'\, \psi.
$$
iii) Every second-order operator ${\cal F}_{\gamma}\otimes {\cal
F}_{\lambda} \longrightarrow {\cal F}_{\gamma+\lambda+2}$ is given
by the composition of the de Rham operator and the Poisson
bracket:
$$
\begin{array}{lcl}
\varphi (dx)^{\gamma} \otimes \psi (dx)^{\lambda} & \mapsto & d \{
\varphi,\psi \}\,dx\quad
\mbox{for } \lambda+\gamma+1=0,\\[2mm]
\varphi (dx)^{\gamma} \otimes \psi(dx)^{\lambda}
 & \mapsto & \{d\,\varphi,\psi\}(dx)^{\lambda+2} \quad
\mbox{for } \gamma=0,\\[2mm]
\varphi(dx)^{\gamma} \otimes \psi (dx)^{\lambda} & \mapsto &
\{\varphi,d\,\psi\}(dx)^{\gamma+2} \quad \mbox{for } \lambda=0.
\end{array}
$$

iv) Every third-order invariant operator ${\cal F}_{\gamma}\otimes
{\cal F}_{\lambda} \longrightarrow {\cal F}_{\gamma+\lambda+3}$ is
given by the composition of the de Rham operator and the Poisson
bracket. Moreover, there exists another operator, now called Grozman
operator, which is not given by this composition. The list is as
follows:
$$
\begin{array}{lcll}
\varphi\otimes \psi(dx)^{-2} & \mapsto &d\,\{d\,\varphi,\psi\}\,dx&
\mbox{for } (\gamma,\lambda)=(0,-2)\\[2mm]
\varphi(dx)^{-2} \otimes \psi& \mapsto & d\,\{\varphi,d\,\psi\}\,dx&
\mbox{for } (\gamma,\lambda)=(-2,0)\\[2mm]
\varphi \otimes \psi& \mapsto &\{d\,\varphi,d\,\psi\}(dx)^{3} &
\mbox{for } (\gamma,\lambda)=(0,0)\\[2mm]
\varphi (dx)^{-2/3} \otimes \psi(dx)^{-2/3}& \mapsto & \left (\left
|
\begin{array}{cc}
\varphi &\varphi'''\\
\psi & \psi'''
\end{array}
\right | -\frac{3}{2} \left |
\begin{array}{cc}
\varphi' &\varphi''\\
\psi' & \psi''
\end{array}
\right |\right ) (dx)^{5/3}
& \mbox{for } (\gamma,\lambda)=(-\frac{2}{3},-\frac{2}{3})\\[2mm]
\end{array}
$$

The list of ${\mathfrak l}$-invariant bilinear differential
operators on weighted densities is a classical result of Gordan
\cite{g}. These operators are called ``transvectants'' (see
\cite{bo, gar}). Their explicit expressions are given by
$$
\begin{array}{cll}
{\cal F}_{\gamma}\otimes {\cal F}_{\lambda}& \longrightarrow
_{\gamma+\lambda+k}\\[3mm]
\varphi (dx)^{\gamma}\otimes \psi (dx)^{\lambda} & \mapsto
&\displaystyle \sum_{i+j=k} (-1)^{i} {k \choose
i}\frac{(2\gamma-i)\cdots (2\gamma -k+1)}{(2\lambda-j)\cdots
(2\lambda -k+1)} \,\varphi^{(i)}\,\psi^{(j)}(dx)^{\gamma+\lambda+k}.
\end{array}
$$
Now, we investigate bilinear differential operators that are
invariant with respect to the Lie subalgebras ${\mathfrak k}_1$
and ${\mathfrak k}_2$.
\begin{thm} \label{cns} Up to a constant, the only ${\mathfrak
k}$-invariant bilinear differential operators on tensor densities
are: (i) the product of tensor densities, (ii) the Poisson bracket, (iii)
the composition of the de Rham operator and Poisson bracket, (iv) the
Grozman operator.
\end{thm}
{\bf Proof.} Any bilinear differential operator ${\cal A}:{\cal
F}_{\gamma}\otimes {\cal F}_{\lambda}\longrightarrow {\cal
F}_{\mu}$ can be expressed as
\begin{equation}
\label{eqq} {\cal A}(\varphi,\psi)=\sum_{i+j\leq k}
c_{i,j}\,\varphi^{(i)}\psi^{(j)}(dx)^\mu,
\end{equation}
where $c_{i,j}$ are smooth functions and the superscript $(i)$
stands for the $i$-th derivative. The ${\mathfrak k}$-invariance of
the operator above reads as follows: for any
$(\varphi(dx)^\gamma,\psi(dx)^\lambda) \in {\cal
F}_\gamma\otimes{\cal F}_\lambda$ and $X\frac{d}{dx} \in\mathfrak
k$, we have
 \begin{equation}
 \label{eq}
L_{X\frac{d}{dx}}^\mu{\cal A}(\varphi,\psi)={\cal
A}(L_{X\frac{d}{dx}}^{\gamma}\varphi,\psi)+{\cal A}
(\varphi,L_{X\frac{d}{dx}}^\lambda\psi).
\end{equation}
The ${\mathfrak k}$-invariance (\ref{eq}) for
$X\frac{d}{dx}=\frac{d}{dx}$ implies that $c_{i,j}'=0.$ Thus,
$c_{i,j}$ are constants.

(a) For $k=0,$ the ${\mathfrak k}$-invariance
 (\ref{eq}) for $X\frac{d}{dx}=\sin(sx) \frac{d}{dx}$ and $\cos(sx)
\frac{d}{dx}$ (or $\sinh(sx) \frac{d}{dx}$ and $\cosh(sx)
\frac{d}{dx}$) is equivalent to the following equation
\[
\left(-\gamma -\lambda +\mu \right) \,c_{0,0}=0.
\]
Therefore, $\mu=\gamma +\lambda$ and the corresponding operator is
the multiplication operator given as in Grozman's list.

 (b) For $k=1,$ the ${\mathfrak k}$-invariance is equivalent to
the following system
\begin{equation}
\label{eqqq}
\begin{array}{lcllcl}
(1+\gamma+\lambda-\mu)\,c_{1,0}&=&0,&(1+\gamma+\lambda-\mu)\,c_{0,1}&=&0,\\
(\gamma+\lambda -\mu)\, c_{0,0}&=&0,& \gamma \,c_{0,1}+\lambda
\,c_{1,0}&=&0.
\end{array}
\end{equation}
If $1+\gamma+\lambda-\mu\not=0,$ then $c_{1,0}=c_{0,1}=0.$ Thus
the operator is $0$-order, but this case has already been studied
in Part (a). Therefore, $1+\gamma+\lambda-\mu=0.$ The third
equation implies that $c_{0,0}=0$ which implies that the operator
(\ref{eqq}) is homogeneous. Now the last equation admits a
solution for all $\lambda$ and $\gamma$. The corresponding
operator coincides (up to a factor) with the Poisson bracket.

(c) For $k=2,$ the ${\mathfrak k}$-invariance is equivalent to the
system (\ref{eqqq}) and the following system:
\[
\begin{array}{rclrclrcl}
\lambda \,c_{1,1}+(1+2\gamma)\,c_{2,0}&=&0,& \gamma\,
c_{1,1}+(1+2\lambda)\,c_{0,2}&=&0,& \gamma\, c_{2,0}+\lambda \,
c_{0,2}&=&0,\\[2mm]
(2 + \gamma + \lambda - \mu)\,c_{1, 1}&=&0,&(2 + \gamma + \lambda -
\mu)\, c_{2, 0}&=& 0,&(2 + \gamma + \lambda -\mu)\, c_{0, 2} &=&0.
\end{array}
\]
If $2 + \gamma + \lambda -\mu\not=0,$ then the constants
$c_{2,0}=c_{0,2}=c_{1,1}=0.$ In this case the operator (\ref{eqq})
becomes $1$-order which has already been studied in Part (b). If
$2 + \gamma + \lambda -\mu=0,$ then the constants
$c_{1,0}=c_{0,1}=c_{0,0}=0$ and the operator becomes homogeneous.
The only values of $\lambda$ and $\gamma$ for which the system is
consistent are given as in Part iii) of Grozman's list.

(d) For $k=3,$ we proceed as before and the operator (\ref{eq}) is
given as in Grozman's list.

(e) Let us prove that for $k\geq 4$ there are no operators
invariant with respect to the Lie algebra ${\mathfrak k}.$ We
proceed by induction; of course, a direct proof can be checked for
$k=4.$ The invariance property implies that the coefficients of
the components $X' \varphi^{(s)}\psi^{(k-s)}$ which should vanish
is given by the system
\[
(k+\gamma+\lambda-\mu)\,c_{s,k-s} \quad \mbox{for} \quad
s=0,\ldots,k.
\]
If $k+\gamma+\lambda-\mu\not=0,$ then all the constants
$c_{s,k-s}=0.$ Thus, the operator (\ref{eqq}) is $(k-1)$-order and
the induction assumption assures the non-existence. Suppose then
that $k+\gamma+\lambda-\mu=0.$ The invariance implies
that the components of $X'\varphi^{(s)}\psi^{(t)}$, where $s+t<k,$
are given by
\[
(s+t+\gamma+\lambda-\mu)\,c_{s,t}=0 \quad \mbox{for} \quad
s+t=0,\ldots,k-1.
\]
Thus, all $c_{s,t}=0$ for all $s+t<k.$ Hence the operator
(\ref{eqq}) becomes homogeneous, namely
\[
\sum_{i+j=k}
c_{i,j}\,\varphi^{(i)}\psi^{(j)}(dx)^{\gamma+\lambda+k}.
\]
Now, the invariance gives a system equivalent to the system coming
from the invariance property with respect to the Lie algebra
$\mathrm{Vect}(M).$ The result follows from Grozman's
classification.\\
\cqfd
\subsection{Proof of Theorem \ref{t1}.}
\label{sw}
\begin{lem}\cite{lo}
Every 1-cocycle $c$ on $\mathrm{Vect}(M)$ with values in
$\mathrm{Hom}_{\mathrm{diff}} ({\cal F}_\lambda, {\cal F}_\mu)$ is
differentiable.
\end{lem}
\begin{lem}
Every 1-cocycle $c$ on $\mathrm{Vect}(M)$ that vanishes on
${\mathfrak k}$ is necessarily ${\mathfrak k}$-invariant.
\end{lem}
{\bf Proof.} The 1-cocycle property reads
$$
L_{X\frac{d}{dx}}\,c(Y,A)-L_{Y\frac{d}{dx}}\,c(X,A)=c([X,Y],A),
$$
for any $X\frac{d}{dx}, Y\frac{d}{dx}\in \mathrm{Vect}(M)$ and $A\in
\mathrm{Hom}_{\mathrm{diff}} ({\cal F}_\lambda, {\cal F}_\mu)$. As
$c$ vanishes on ${\mathfrak k},$ then
$$
L_{X\frac{d}{dx}}\,c(Y,A)=c([X,Y],A)\quad \mbox{ for any } X\frac{d}{dx}\in {\mathfrak k}.
$$ This is just the $\mathfrak k$-invariance property.\\
\cqf

According to Theorem \ref{cns}, there are only two ${\mathfrak
k}$-invariant differential operators
$$
\mathrm{Vect}(M)\otimes {\cal F}_{\lambda}\longrightarrow {\cal
F}_{\mu},
$$ and they are given by the Poisson bracket $\{\cdot,
\cdot\}$ for $\gamma=-1$ and the operator $\{\cdot,d\} $ for
$(\gamma,\lambda)=(-1,0).$ Let us describe under what conditions
the operator $\{\cdot,d\}$ is a $1$-cocycle.
For any $\psi \in {\cal F}_{0}$ and $X\frac{d}{dx}, Y\frac{d}{dx}\in
\mathrm{Vect}(M)$, we have
$$
\{[X,Y],d\,\psi\}=L_{X\frac{d}{dx}}\{Y,d\,\psi\}-\{Y,d\,L_{X\frac{d}{dx}}\psi\}-
L_{Y\frac{d}{dx}} \{X,d\,\psi\} +\{X,d\,L_{Y\frac{d}{dx}}\psi\}.
$$
The cocycle condition implies that the coefficients $c_{2,0},
c_{0,2}$ and $c_{1,1}$ of the operator $\{\cdot,d\}$ satisfy the
system
\begin{eqnarray*}
(\mu-1)\,c_{2,0}=0,\quad \,c_{0,2}=0,\quad (2-\mu)\,c_{0,2}=0.
\end{eqnarray*}
Therefore, the operator $\{\cdot ,d\}$ is never a 1-cocycle. The
same holds for the Poisson bracket $\{\cdot, \cdot\}$, and
therefore we have just proved that
$\mathrm{H}^{1}(\mathrm{Vect}(M),{\mathfrak k};
\mathrm{Hom}_{\mathrm{diff}}({\cal F}_\lambda,{\cal F}_\mu))= 0$
and (i) is proven.

To prove ii), let $c$ be a 1-cocycle on $\mathrm{Vect}(M)$
with value in $\mathrm{Hom}_{\mathrm{diff}}({\cal F}_\lambda,{\cal
F}_\mu), $ and denote by $\widehat{c}$ its restriction to
${\mathfrak k}$. We will prove that the map $c\mapsto \widehat{c}$
is an isomorphism. Let $c_1$ and $c_2$ be two 1-cocycles such that
$\widehat{c}_1=\widehat{c}_2$. It follows that the difference
$c_1-c_2$, which is also a 1-cocycle, vanishes on ${\mathfrak k}$.
Part i) of the Theorem ensures that it is a coboundary, namely
$$
c_1(X)-c_2(X)=L_{X\frac{d}{dx}}\circ A-A\circ L_{X\frac{d}{dx}}
\quad \mbox{for all}\,\, X\frac{d}{dx}\in
\mathrm{Vect}(M),
$$
where $A$ is a certain operator in
$\mathrm{Hom}_{\mathrm{diff}}({\cal F}_\lambda,{\cal F}_\mu)$.
Thus, $c_1\equiv c_2$.

Let $s$ be a 1-cocycle on ${\mathfrak k}$ with values in
$\mathrm{Hom}_{\mathrm{diff}}({\cal F}_\lambda,{\cal F}_\mu)$. We
will prove that there exists a 1-cocycle on $\mathrm{Vect}(M)$, say
$c,$ such that $\widehat{c}=s$. To do that, we write the 1-cocycle
$s$ as follows. For any $X\frac{d}{dx}\in \mathrm{Vect}(M)$ and
$\psi (dx)^{\lambda}\in {\cal F}_\lambda,$ we put (where $c_{i,j}$
are smooth functions)
$$
s(X,\psi)=\sum_{i,j}c_{i,j}\,X^{(i)}\,\psi^{(j)}(dx)^{\mu}.
$$
The cocycle property for $X\frac{d}{dx}=\cos(sx)\frac{d}{dx}$ (or
$\cosh(sx)\frac{d}{dx})$ is equivalent to the following system
\small{
\begin{eqnarray*}
i(i+j-1-\mu+\lambda)\,c_{i,j}-(j+\lambda i)
\left(^{j+i-1}_{\,\,\,\,\,\,\,i} \right )
c_{1,j-1+i}&=&0\quad \mbox{for } i+j\leq k\mbox{ and } i>1,\\
 \frac{j-i}{j}
 \left(^{j+i-1}_{\,\,\,\,\,\,\,i}
\right )c_{i+j-1,l}+\left (
 \left(^{l+i-1}_{\,\,\,\,\,\,\,i}
\right ) + \lambda \left(^{l+i-1}_{\,\,\,i-1\,\,} \right ) \right
)c_{j,l-1+i}&&\\-\left (
 \left(^{l+j-1}_{\,\,\,\,\,\,\,j}
\right ) + \lambda \left(^{l+j-1}_{\,\,\,j-1\,\,} \right ) \right
)c_{i,l-1+j}&=&0\quad \mbox{for }
 l+i+j\leq k+1\mbox{ and } j>i>1,\\
(\mu-i-\lambda)\,c_{0,i}-c_{1,i}'&=&0 \quad \mbox{for } 0<i\leq k,\\
(j-1-\mu+l+\lambda) \,c_{j,l}-\left (
 \left(^{l+j-1}_{\,\,\,\,\,\,\,j}
\right ) + \lambda \left(^{l+j-1}_{\,\,\,j-1\,\,} \right ) \right
)c_{1,l-1+j}&=&0 \quad \mbox{for } j+l\leq k.
\end{eqnarray*}
}
Now, a direct computation proves that the operator defined by
$\sum_{i,j}c_{i,j}\,X^{(i)}\,\psi^{(j)}(dx)^{\mu},$ where $c_{i,j}$
are as above, satisfies the 1-cocycle property on the entire
$\mathrm{Vect}(M)$. In other words, any 1-cocycle on ${\mathfrak k}$
is certainly a 1-cocycle on $\mathrm{Vect}(M)$.

The proof of (ii) for $M=\mathbb{R}$ follows from computations of
$\mathrm{H}^1(\mathrm{Vect}(\mathbb{R}),
\mathrm{Hom}_{\mathrm{diff}}({\cal F}_\lambda,{\cal F}_\mu))$ due
to Feigin-Fuchs \cite{f}. Our Theorem is proven.

\section{Quantization equivariant with respect to
${\mathfrak k}$}
\quad The modules of linear differential operators on the spaces of
tensor densities on a smooth manifold are classical objects (see
\cite{w}). Recently, they have been intensively studied in a series
of papers (see \cite{b2, bg, bo, do, do2, gar, lmt, lo}). Let ${\cal
D}_{\lambda,\mu}^k$ be the space of order $k$ differential operators
$A:{\cal F}_{\lambda}\to{\cal F}_{\mu}$ endowed with the
$\mathrm{Vect}(M)$-module structure by the formula
\begin{equation}
\label{act} L_{X\frac{d}{dx}}^{\lambda,\mu}(A)=L_{X\frac{d}{dx}}^\mu
\circ A-A\circ L_{X\frac{d}{dx}}^\lambda.
\end{equation}

Let $\delta=\mu-\lambda$. For first-order differential operators,
 it was proven in \cite{gar} that there exists a map
$$
{\cal F}_{\delta-1}\longrightarrow {\cal D}_{\lambda, \mu}^{1}\qquad
a_{1}(dx)^{\delta-1} \mapsto  \left (a_{1}\frac{d}{dx}
+\frac{\lambda}{1-\delta}a_{1}'\right ) dx^{\delta}
$$
equivariant with respect to $\mathrm{Vect}(M)$. However, for
second-order differential operators, it was proven in \cite{cmz,
gar} that there exists a map
$$
{\cal F}_{\delta-2}  \longrightarrow  {\cal D}_{\lambda,
\mu}^{2}\qquad a_2(dx)^{\delta-2} \mapsto \left (
a_{2}\frac{d^2}{dx^2}+\frac{1+2\lambda}{2-\delta}a_2'\frac{d}{dx}
+\frac{\lambda(1+2\lambda)}{(2-\delta)(3-2\delta)}a_{2}''\right
)dx^{\delta}
$$
equivariant with respect to ${\mathfrak l}_0$ but there is no map
${\cal F}_{\delta-2} \longrightarrow  {\cal D}_{\lambda, \mu}^{2}$
equivariant with respect to the whole Lie algebra
$\mathrm{Vect}(M)$. For ${\mathfrak k},$ we have the following
\begin{thm}
A differential operator
$$
Q:{\cal F}_{\delta-2}\oplus {\cal F}_{1-\delta}\oplus {\cal
F}_{\delta}\longrightarrow {\cal D}_{\lambda, \mu}^{2}
$$
equivariant with respect to ${\mathfrak k}$ exists only for either
$\mu=1$ and $\lambda\not =-1,$ or $\mu\not =2$ and $\lambda=0.$
\end{thm}
{\bf Proof.} Assume the contrary: Let the map $Q$ exist. We have,
therefore, a diagram
$$
\begin{CD}
{\cal F}_{\delta-2}\oplus {\cal F}_{\delta-1}\oplus {\cal
F}_{\delta}@>Q>> {\cal D}_{\lambda, \mu}^2\\
@VVV @VVL_{X\frac{d}{dx}}^{\lambda, \mu}V\\
{\cal F}_{\delta-2}\oplus {\cal F}_{\delta-1}\oplus {\cal
F}_{\delta}@>Q>> {\cal D}_{\lambda, \mu}^2
\end{CD}
$$
that by our assumption commutes once the action is restricted to
${\mathfrak k}$. Hence the map
$$
Q^{-1}\circ L_{X\frac{d}{dx}}^{\lambda, \mu} \circ Q,
$$
defines a new action on the space ${\cal F}_{\delta-2}\oplus {\cal
F}_{\delta-1}\oplus {\cal F}_{\delta}$. Thus, there exists a map
$c$ such that
$$
Q^{-1}\circ L_{X\frac{d}{dx}}^{\lambda, \mu} \circ
Q=(L_{X\frac{d}{dx}}^{\delta-2},L_{X\frac{d}{dx}}^{\delta-1},
L_{X\frac{d}{dx}}^\delta)+c(X),
$$
where $L_{X\frac{d}{dx}}^\delta$ is the action (\ref{dens}). It is
easy to see that $c$ defines a 1-cocycle on $\mathrm{Vect}(M)$
with values in $\mathrm{End}_{\mathrm{diff}}({\cal
F}_{\delta-2}\oplus {\cal F}_{\delta-1}\oplus {\cal F}_{\delta}),$
vanishing on ${\mathfrak k}$. By Theorem \ref{cns}, the map $c$ is
identically zero, not just a coboundary. Thus, the map $Q$ is not
only $\mathfrak k$-equivariant but also
$\mathrm{Vect}(M)$-equivariant. But by \cite{gar} there is no
quantization map equivariant with respect to the whole Lie algebra
$\Vect(M)$ except for $\mu=1$ and $\lambda\not =-1,$ or $\mu\not
=2$ and $\lambda=0.$\\
\cqfd
\subsection*{Conclusion}
\quad Even though the Lie subalgebras ${\mathfrak k}_1$ and
${\mathfrak k}_2$ are of finite dimension, the (differential)
equivariant quantization map does not exist. We have shown in this
paper that the existence is related to the nature of the vector
fields that generate a Lie subalgebra, rather than the number of
generators. We believe that the dimension of a Lie subalgebra will
determine only the uniqueness of the quantization map -- whenever
it exists -- as already pointed out in \cite{dlo}.

\subsection*{Acknowledgments}
We are grateful to Profs. L. Ayari, D. Leites, V. Ovsienko, and
the referees for their suggestions and remarks.

\label{lastpage}

\end{document}